\documentstyle[12pt]{article} 

\def\basicalgebra{{\cal A}} 
\def\subalgebra{{\cal B}} 
 
\def\constants{{\bf C}} 
\def\difforms{{\bf\Omega}} 
\def\conjtoder{{\bf W}} 
\def\derivatives{{\cal V}} 
\def\der{{\sf Der}}
\def\annihilator{{\sf Ann}}
\def\image{{\sf Im}}
\def\kernel{{\sf Ker}}
\def\inclusion{{\sf in}}

\def\scalars{{\cal S}} 
\def\field{{\bf F}} 
\newcommand{\dff}\sc  
 
\def\be{\begin{equation}\label} 
\def\ee{\end{equation}} 
\def\zentrum{{\cal Z}} 
\def\jpi{{j_\pi}} 

\def\proof{\paragraph{Proof }} 
\def\endproof{\hfill$\Box$\medskip}

\newtheorem{proposition}{Proposition} 
\newtheorem{definition}{Definition}

\title{DUAL STRUCTURES IN NON-COMMUTATIVE DIFFERENTIAL ALGEBRAS} 

\author{G.N.Parfionov\thanks{Department of Mathematics, SPb UEF,
Griboyedova 30/32, 191023, St-Petersburg, Russia (address for
correspondence)}, R.R.Zapatrin\thanks{Division of Mathematics,
Instituto per la Ricerca di Base, I-86075, Monteroduni, Molise, Italy}} 

\date{} 

\begin{document} 
\maketitle 

\begin{abstract} The non-commutative algebraic analog of the moduli
of vector and covector fields is built. The structure of
moduli of derivations of non-commutative algebras are studied. 
The canonical coupling is introduced and the conditions
for appropriate moduli to be reflexive are obtained. 
\end{abstract} 

\section*{FOREWORD} 

The duality problem we are going to tackle stems from the 
non-commutative generalization of differential geometry of 
manifolds. The variety of structures it deals with is defined in 
terms of two basic objects: the algebra $\basicalgebra$ of smooth 
functions on a manifold and the Lie algebra $\derivatives$ of 
smooth vector fields. In standard differential geometry 
$\derivatives$ is always a {\em reflexive} $\basicalgebra$-module 
which enables the tensor analysis to be successively built. 
Whereas it is worthy to note that the fundamental geometrical 
notions are formulated in pure algebraic terms of {\em 
commutative} algebra $\basicalgebra$ \cite{r2}. 
The attempts of direct generalization of these notions starting 
from a non-commutative algebra $\basicalgebra$ produce non-trivial 
algebraic problems: we focus on two of them. 

The first is to introduce an analog of the module of vector
fields: in classical geometry it is the $\basicalgebra$-module
$\derivatives = \der \basicalgebra$, meanwhile when
$\basicalgebra$ is non-commutative $\der \basicalgebra$ is not
$\basicalgebra$-module.  To overcome it, we have to restrict
$\der \basicalgebra$ to some its Lie subalgebra $\derivatives$,
and introduce the set of scalars $\scalars$ {\em with respect
to} $\derivatives$. In this paper we shall use the center
$\zentrum$ of $\basicalgebra$ as the set of scalars. 

The second problem is to introduce an appropriate definition of
the dual module of $\derivatives$. The matter is that the
standard dualization borrowed from the theory of moduli turns
out to be incompatible with the generalization of forthcoming
geometrical notions such as, for instance, Cartan differential.
The guideline to solve the second problem is to choose the dual
module so that the obtained pair of dual objects would be
reflexive. The standard definitions related with rings and
moduli are borrowed from \cite{r1}. 

The paper is structured as follows. In section \ref{s1} the
notion of differential algebra is introduced as the couple
(associative~algebra, module of~vector~fields). In section
\ref{s2} different ways to build the dual to the module of
'vector fields' are discussed and the choice of appropraite
definition is motivated. In section \ref{s3} the conditions for
the module of vector fields to be refleive are introduced. 

\section{DIFFERENTIAL ALGEBRAS}\label{s1}

Let $\field$ be a field with zero characteristic, and consider a 
(non-commutative, in general) associative algebra with the unit 
element over $\field$. Denote by $\der \basicalgebra$ the set of 
the set of {\dff derivatives} of $\basicalgebra$, that is the 
$\field$-linear mappings $v:\basicalgebra\longrightarrow \basicalgebra$ which 
enjoy the Leibniz rule: for any $a,b\in \basicalgebra$ 
\[ 
v(ab) = va\cdot b + a\cdot vb 
\] 

$\der \basicalgebra$ is the Lie algebra over $\field$ with
respect to the bracket operation $[u,v]a = uv(a)-(vu)a$. The
action of $\der \basicalgebra$ on $\basicalgebra$ induces the
Galois connection between subsets of $\basicalgebra$ and $\der
\basicalgebra$ defined as follows. For each $\subalgebra \subseteq 
\basicalgebra$ ($\derivatives\subseteq \der \basicalgebra$, resp.) 
define the polar $\subalgebra^c$ ($\derivatives^c$, resp.) as:
\be{f2} 
\subalgebra^c =\{ v\in \derivatives\mid  \forall b\in \subalgebra 
\quad vb=0\} \qquad \derivatives^c =\{ a\in \basicalgebra\mid  \forall 
v\in \derivatives \quad va=0\} 
\ee 

It follows from the general theory of Galois connections that 
$\derivatives^{ccc} = \derivatives^c$ and $\subalgebra^{ccc} = 
\subalgebra^c$. 

\begin{proposition}\label{l11} \begin{enumerate} 
\item \label{i111} For any $\derivatives\subseteq \der 
\basicalgebra$ the set $\derivatives^c$ is the subalgebra of 
$\basicalgebra$.
\item \label{i112} For any $\subalgebra\subseteq \basicalgebra$ 
the set $\subalgebra^c$ is the Lie subalgebra of $\der 
\basicalgebra$.
\end{enumerate} 
\end{proposition} 

\proof is yielded by direct verification of appropriate 
definitions. \endproof

     For any $a\in \basicalgebra$, $v\in \derivatives$ the linear 
mapping $av:\basicalgebra\longrightarrow \basicalgebra$ is defined as 
\be{f12} 
(av)b=a\cdot vb
\ee

\begin{definition}\label{d11} A {\dff differential algebra} is a pair 
$(\basicalgebra,\derivatives)$ where $\basicalgebra$  is  an 
algebra and $\derivatives\subseteq \der \basicalgebra$ are such that 
$\derivatives=\derivatives^{cc}$. 
\end{definition} 

It follows immediately from Proposition \ref{l11}  
that $\derivatives$ is always the Lie subalgebra of $\der 
\basicalgebra$. There are two algebras associated with 
$(\basicalgebra,\derivatives)$:
\be{f3}  
\begin{array}{rclcl} 
\constants &=& \derivatives^c &\qquad& \hbox{ the algebra of 
constants} 
\end{array} 
\ee 

It turns out that a differential algebra can be equivalently
defined as the pair $(\basicalgebra,\constants)$, with
$\constants =\constants^{cc}\subseteq \basicalgebra$, which
coincides with the above definition \ref{d11} by putting
$\derivatives=\constants^c$.  In any differential algebra
$\derivatives$ possesses the natural structure of the module
over the commutative algebra $\zentrum$ (the center of
$\basicalgebra$).

\section{DUALITY}\label{s2}  

There is the {\em canonical} construction of conjugated to a left 
(right) $\zentrum$-module which is the right (resp., left) 
$\zentrum$-module. So, the conjugated to the left $\zentrum$-module 
$\derivatives$ will be the right $\zentrum$-module 
\[ 
\derivatives^\ast = \hom_{\zentrum} (\derivatives,\zentrum) 
\] 

Likewise the  left  $\zentrum$-module $\derivatives^{\ast\ast}$ is 
introduced. There is the canonical homomorphism 
$\derivatives\longrightarrow \derivatives^{\ast\ast}$ of left 
$\zentrum$-moduli which is not isomorphism in general. The reflexive 
moduli (with $\derivatives=\derivatives^{\ast\ast}$) are normally 
"good" objects to deal with. In our paper we are going to fetch 
the reflexivity to the construction by appropriate choice of the 
definition of conjugation taking into account the features of 
$\derivatives\subseteq \der \basicalgebra$. Consider the set:
\[ 
\derivatives^+   = \hom_{\zentrum} (\derivatives,\basicalgebra)
\] 

\noindent of all \zentrum-linear forms on $\derivatives$ taking values in 
$\basicalgebra$ (rather than in $\zentrum$). Call the elements of 
$\derivatives$ vectors and the elements of $\derivatives^+$ 
covectors. The following propositions show the relevance of such
definition.

\begin{proposition}\label{l21} $\derivatives^+$ possesses the 
natural structure of $\basicalgebra$-bimodule.
\end{proposition} 

\proof The additive structure on $\derivatives^+$ is introduced in 
the standard way, and the right action of the elements of 
$\basicalgebra$ is defined as: 
\[ 
(\omega \cdot a)(v) =  \omega (v)a \qquad \qquad  \omega \in 
\derivatives^+\quad\hbox{,}\quad a\in A, v\in V 
\] 

\noindent making $\derivatives^+$ right $\basicalgebra$-module. 
For any $a\in \basicalgebra$, $\omega \in \derivatives$ consider 
the $\basicalgebra$-valued function $(a \omega )(v) = a \omega (v)$ 
on $\derivatives$. To check that $a \omega \in \derivatives^+$, it 
suffices to check its $\zentrum$-linearity. For each $s\in 
\zentrum$ consider the discrepancy $\delta = a \omega (sv) - sa 
\omega (v) = [a,s] \omega (v) =0$, since $s\in\zentrum$. 
Finally, for each $v\in \derivatives$ 
we have $(a \omega (v))b = a \omega (v)b = a( \omega (v)b)$, hence
\[ 
(a \omega )b = a( \omega b) \qquad\qquad  a,b \in  
\basicalgebra\quad\hbox{,}\quad \omega  \in  \derivatives^+ 
\] 

\noindent therefore $\derivatives^+$   is 
$\basicalgebra$-bimodule. \endproof 

\paragraph{Remarks.} (i). This unexpected feature of 
$\derivatives^+$ was enabled by that 
$\derivatives$ is the module of derivations. 

\noindent (ii). In general $a \omega \neq \omega a$. 
\medskip

Each element $a\in \basicalgebra$ canonically induces the 
$\field$-linear mapping $da:\derivatives\longrightarrow 
\basicalgebra$ defined as $da(v) = va$. The next statement
solves the problem of generalization of Cartan's differentials. 

\begin{proposition}\label{l22} \begin{enumerate} 
\item \label{i221} $da\in \derivatives^+$ for any $a\in 
\basicalgebra$. 
\item \label{i222} \( d(ab) = da\cdot b+a\cdot db \)  
\item \label{i223} \( \ker d\quad =\quad \constants\quad = 
\quad \derivatives^c \)  
\end{enumerate} 
\end{proposition} 

\proof \ref{i221} $da(sv) = sva = sda(v)$, \ref{i222} is verified 
by direct checking, \ref{i223} follows directly from the
definitions (\ref{f2},\ref{f3}). \endproof 
\medskip

We shall call the elements of $\derivatives^+$ of the form $da$ 
{\dff differentials}. Consider the $\basicalgebra$-bisubmodule of 
$\derivatives$ generated by differentials: 
\[ \difforms = \{  \omega \in \derivatives^+\mid   \omega = \sum_i\rho_i 
da_i q_i\} \] 

\noindent where $i$ ranges over a finite set and call $\difforms$ 
the {\dff module of differential forms}. Introduce the set 
\be{f4}  
\conjtoder = 
\hom_{\basicalgebra,\basicalgebra}(\derivatives^+,\basicalgebra)
\ee 

\noindent of all $\basicalgebra$-bilinear homomorphisms
$\derivatives^+ \longrightarrow \basicalgebra$. 

\begin{proposition}\label{l24} $\conjtoder$ is $\zentrum$-bimodule. 
\end{proposition} 

\proof For $w\in \conjtoder$, $s\in \zentrum$ define the mapping 
$sw:\derivatives^+ \longrightarrow \basicalgebra$ as $(sw) \omega 
=s\cdot w \omega$. In the similar way the right action of
$\zentrum$ is introduced. Then prove that $swt$ is 
$\basicalgebra$-bilinear for any $s,t \in \zentrum$: 
\[ 
(swt)(a \omega b) = s\cdot a\cdot w \omega \cdot b = 
a\cdot (s w \omega t)\cdot b
\] 
\endproof 

The bilinear form $<\cdot ,\cdot>$ on $\derivatives \times 
\derivatives^+$ with the values in $\basicalgebra$ naturally arises:
\[ 
<v, \omega > =  \omega (v)
\] 
\noindent Note that no confusion occurs: $<sv,a \omega b> = sa<v, 
\omega b> = as<v, \omega >b$ for $s\in \zentrum$. This bilinear form 
enables the canonical homomorphism $j:\derivatives \longrightarrow 
\conjtoder$ defined as: 
\be{f21} 
jv = <v,\cdot > 
\ee 

In the category of moduli there always exists the natural 
homomorphism $\derivatives \longrightarrow 
\derivatives^{\ast\ast}$, which is in general neither isomorphism, 
nor even embedding. In our scheme the r\^ole of 
$\derivatives^{\ast\ast}$ is played by $\conjtoder$, and the 
following holds:

\begin{proposition}\label{l25} The homomorphism $j$ embeds 
$\derivatives$ into $\conjtoder$ as direct summand. 
\end{proposition} 

\proof For each $w\in \conjtoder$ define the mapping 
$\hat{w}:\basicalgebra \longrightarrow \basicalgebra$ as follows: 
$\hat{w}(a) = w(da)$. Then $\hat{w}(ab) = w(da\cdot b)+w(a\cdot db) 
= wa\cdot b+a\cdot \hat{w}b$ , hence $\hat{w}\in \der 
\basicalgebra$. Since $\hat{w}(k)=w(dk)=0$ for any $k\in 
\constants$, $w\in 
\derivatives$ (it follows from the definition of differential 
algebra that $\derivatives=\constants^c$). Denote 
$w\mapsto\hat{w}$ by $\pi$: 
\[ 
\pi (w) = \hat{w} 
\] 
\noindent Evidently, $\pi$ is $\zentrum$-linear, therefore it is the 
homomorphism $\conjtoder \longrightarrow \derivatives$. Then 
consider $(\pi jv)(a) = <v,da> = va$, hence $\pi 
j=1_{\derivatives}$, thus $j$ is the embedding. \endproof 

\begin{proposition} Consider two submoduli of $\conjtoder$: $\image j$ (being
isomorphic to $\derivatives$) and \( N = \kernel \pi \). Then 
\be{f2277} 
\conjtoder  = \derivatives\oplus N
\ee 
\end{proposition} 
 
\proof Consider the endomorphism $\jpi =j\circ\pi$ of $\conjtoder$. Then 
$\kernel \jpi = \kernel \pi = N$ (since $\kernel j=0$), and 
$\image \jpi  = \image j$ (since $\image \pi =\derivatives$). 
Finally, $\jpi^2 = \pi j\pi j = \jpi$, hence $\conjtoder = \image 
j\oplus N$. \endproof 

Since $j$ is the embedding, in the sequel we shall identify 
$\derivatives$ with its image, so

\paragraph{Corollary.} $N$ is exhausted by annihilators of 
differential forms in $\derivatives^+$: 
\be{f23} 
N = \annihilator  \difforms  = \{ n\in  \conjtoder \mid  \forall  
\omega \in  \difforms \quad \omega (n)=0\}
\ee

\section{REGULARITY AND REFLEXIVITY} \label{s3} 

Consider an important special case.

\begin{proposition}\label{l26} If $\derivatives$ is the finitely 
generated free module, then $\derivatives$ is canonically 
isomorphic to $\conjtoder$. 
\end{proposition} 

\proof Recall that $\zentrum$ is a subalgebra of 
$\basicalgebra$, hence 
\[ \derivatives^\ast  = \hom_s (\derivatives,\scalars) \subseteq  
\derivatives^+ 
\] 
\noindent Let $v_1,\ldots, v_n$ be a basis of $\derivatives$. 
Build its dual basis $\omega^1,\ldots,\omega^n$ in 
$\derivatives^\ast$: $\omega^i(v_k)=\delta^i_k$. Note that 
$\omega^i (v)\in \scalars$ and prove that $\{\omega^i\}$ is the 
basis of $\derivatives^+$. For any $\omega \in \derivatives^+$ we 
have $\omega (v) =  \omega^i (\sum \omega (v)\cdot v_i) = 
\sum\omega^i (v)\cdot  \omega(v_i)$, hence 
\be{f24} 
\omega^i  = \sum\omega \cdot \omega (v_i)
\ee 

\noindent Moreover, for any $a\in \basicalgebra$, $a \omega^i(v_k) = 
a\delta^i_k = \delta^i_k a = ( \omega^i a)(v_k )$, hence 
\be{f25} 
a \omega^i = \omega^i a
\ee 
\noindent therefore $\forall w\in \conjtoder$ we have 
$aw( \omega^i )=w(\omega^i )a$, and thus $w( \omega^i )\in \zentrum 
(\basicalgebra)\subseteq \scalars$. 

\noindent Denote $w_i = j(v_i)$, and show that any $w\in  
\conjtoder$  can be decomposed over $w_i$: for any $\omega \in 
\derivatives^+$ \ $w( \omega ) = \sum w( \omega^i) \omega (v_i) =$ 
\[ 
\sum \omega (w( \omega^i)v_i) =  
\omega (\sum w( \omega^i)v_i) = 
<\sum w(\omega^i)v_i, \omega> = 
\] 
\[
\sum w( \omega^i)<v_i, \omega > = 
\sum w(\omega^i)\cdot jv_i(\omega) = 
\sum w( \omega^i)\cdot w_i(\omega) 
\] 

\noindent So, $\{w \}$  is the generating set for  $\conjtoder$, 
hence  $\conjtoder$  coincides with the image of the injection 
$j$, which completes the proof. \endproof

Since  $\conjtoder$  splits into the direct sum (\ref{f2277}), 
we can also split out the module $\derivatives^+$. 
Namely, consider the set $R$ of annihilators of $N$: 
\[ 
R = \annihilator N = \{  \omega \in \derivatives^+ \mid  \forall 
n\in N\quad n( \omega )=0\} 
\] 

\noindent and call them {\dff regular covectors}. $R$ is the 
double annihilator of the module of differential forms: 
\[ 
R = \annihilator (\annihilator  \difforms ) \subseteq  
\derivatives^+ 
\] 

A differential algebra $(\basicalgebra,\derivatives)$ is said to 
be {\dff reflexive} if its set of covectors is closed: 
$\derivatives^+ =R$. The adequacy of using the term 'reflexive' is
corroborated by the following 

\begin{proposition}\label{l31} The following two conditions are 
equivalent: \begin{enumerate} 
\item \label{i311} $(\basicalgebra, \derivatives)$ is reflexive. 
\item \label{i312} $\conjtoder =\derivatives$. 
\end{enumerate} 
\end{proposition} 

\proof is straightforward. \endproof 

\paragraph{Remark.} It follows from  Proposition  \ref{l26} that 
any differential algebra with a free finitely generated module 
of vectors is always reflexive. However, the converse is not true: 
the counterexample is yielded by the algebra $\basicalgebra$ of 
smooth functions on the 2-dimensional sphere with 
$\derivatives=\der \basicalgebra$ (since each covector is 
differential form). 
\medskip 

Denote by $R^+$ the dual to $R$: 
\be{f6}  
R^+ = \hom_{\basicalgebra,\basicalgebra} (R,\basicalgebra) 
\ee 

Since $R\subseteq \derivatives^+$ and $\conjtoder 
=\hom_{\basicalgebra,\basicalgebra} (\derivatives^+ 
,\basicalgebra)$, there is the natural homomorphism $\beta:\conjtoder  
\longrightarrow R^+$ such that: 
\[ \beta(w) = w\mid_{\displaystyle R} \] 
for any $w\in \conjtoder$ (called the restriction
homomorphism). 

\begin{proposition}\label{l32} $\kernel \beta  = N$. 
\end{proposition} 

\proof $\kernel \beta  = \{ w\in  \conjtoder \mid \forall  \omega 
\in R w( \omega )=0\}  = \annihilator R = \annihilator 
(\annihilator (\annihilator  \difforms )) = \annihilator  
\difforms  = N$ by virtue of (\ref{f23}). \endproof 

\paragraph{Corollaries.}  (i).  The  homomorphism   $\beta$  can 
be canonically decomposed into canonical projection $\jpi$  and 
monomorphism $i$: 
\be{3.1}
\unitlength=1mm
\begin{picture}(25.00,18.00)
\put(15.00,0.00){\makebox(0,0)[cc]{$\conjtoder / N$}}
\put(0.00,15.00){\makebox(0,0)[cc]{$\conjtoder$}}
\put(25.00,15.00){\makebox(0,0)[cc]{$R^+$}}
\put(2.00,15.00){\vector(1,0){20.00}}
\put(2.00,12.00){\vector(1,-1){10.00}}
\put(18.00,2.00){\vector(1,3){3.33}}
\put(3.00,4.00){\makebox(0,0)[cc]{$\rho$}}
\put(22.00,6.00){\makebox(0,0)[cc]{$i$}}
\put(10.00,18.00){\makebox(0,0)[cc]{$\beta$}}
\end{picture}
\ee 

\noindent (ii). $\derivatives^+ =R$ whenever $\beta$  is 
epimorphism. 

\noindent (iii). $\derivatives$ is the submodule  of  $R^+$  
(since  $\conjtoder /N=\derivatives$). Moreover, the following 
holds:

\begin{proposition}\label{l33} $\derivatives$ is the direct 
summand of $R^+$. 
\end{proposition} 

\proof We already have the injection $i:\derivatives  
\longrightarrow R^+$ (\ref{3.1}). To complete the proof, it 
suffices to build 
the projection $\pi :R^+\longrightarrow \derivatives$ such that 
$\pi i=id_\derivatives$. Define for $x\in R^+$ the mapping $\pi 
x:\basicalgebra \longrightarrow \basicalgebra$ 
\[ 
\pi x(a) = x(da)
\] 

\noindent Like in the Proposition \ref{l25} it  can  be  proved  that 
$\pi x\in \derivatives$.  Now 
for $v\in \derivatives$ consider $\pi (iv)(a) = iv(da) = \beta jv(da) = 
jv(da)$ since $\beta$  is the 
restriction homomorphism, hence $\pi (iv)(a) =  jv(da)  =  v(a)$  by  the 
definition of the mapping $j:\derivatives \longrightarrow  
\conjtoder$  (\ref{f21}). \endproof 

\begin{proposition}\label{34} If $\basicalgebra$ considered 
bimodule over $\basicalgebra$ itself is injective then 
$\derivatives^+=R$. 
\end{proposition} 

\proof It suffices to prove that the injection  $i$ is surjective. 
$\basicalgebra$ is injective, therefore the exactness of the 
sequence of the $\basicalgebra$-bimoduli 
\[ 
0 \longrightarrow R \stackrel{\rm in}{\longrightarrow} 
\derivatives^+
\] 
implies the exactness of 
\[ \hom_{\basicalgebra,\basicalgebra}
(\derivatives^+,\basicalgebra) \longrightarrow 
\hom_{\basicalgebra,\basicalgebra}(R,\basicalgebra) 
\longrightarrow 0
\] 

The first term of the latter sequence is $\conjtoder$
(\ref{f4}), and the
second is $R^+$ (\ref{f6}). The homomorphism between them induced by the
inclusion $\inclusion:R
\longrightarrow \derivatives^+$ is just the restriction 
homomorphism $\beta$, so the sequence 
\[ 
\conjtoder \stackrel{\beta}{\longrightarrow}  R^+ \longrightarrow  
0
\] 
is exact which means that $\beta$  is surjection. 
Finally, the commutativity of the diagram (\ref{3.1}) makes $i$ 
be the surjection. \endproof

\section*{AFTERWORD} 

Two problems concerning the non-commutative generalization of 
differential structure were put forward in the Foreword. 

The first was to suggest a suitable choice of subalgebras 
$\derivatives$ of $\der \basicalgebra$ playing the role of vector 
fields. We did it by introducing the closed subalgebras with 
respect to constants and scalars (section \ref{s1}). 

     Section \ref{s2} prepares the necessary technical tools to 
form the pair of dual objects. Several candidates for being dual 
objects are suggested and some results concerning their structure 
(such as direct decomposability) are obtained. 

     In section \ref{s3} the notion of reflexive differential 
algebra $(\basicalgebra,\derivatives)$ is introduced for which 
the following condition holds: 
\[ 
\derivatives=\hom (\hom 
(\derivatives,\basicalgebra),\basicalgebra)
\] 
where $\hom$ is the set of the morphisms in the aproprate category.

The following differential algebras
$(\basicalgebra,\derivatives)$ were proved to be reflexive
whenever {\em at least one} of 
the following conditions holds: 
\begin{itemize} 
\item $\derivatives$ is finitely generated free module 
\item $\basicalgebra$ is injective $\basicalgebra$-bimodule 
\end{itemize} 

The obtained results are applied to the problem of quantization of 
gravity \cite{r3}. The algebraic account of this problem was first 
done by Geroch \cite{r2} where the commutative case was 
considered. From the general algebraic point of view the idea to
vary the definition of dual object looks rather prospective. In
particular, a representation theory for general partially
ordered sets was built along these lines \cite{d4p,cr}.


\begin{thebibliography}{99} 

\bibitem{r1}  Faith,  C.,  Algebra:  Rings,  Modules  and  Categories,   I, 
Springer, 1973
Faith,  C.,  Algebra II: Ring Theory, Springer, 1976

\bibitem{r2} Geroch, R., 
Einstein Algebras, 
Communications in  Mathematical Physics, 
{\bf 26}, 
271, 
1972 

\bibitem{r3} Parfionov, G.N., Zapatrin, R.R., 
{\it Pointless Spaces  in  General Relativity}, 
International   Journal of   Theoretical   Physics, 
{\bf 34}, 775, 1995

\bibitem{cr} Zapatrin, R.R., 
{\it Les espaces duaux pour les ensembles ordonn\'es\/}, Comptes
   Rendus de l'Acad\'emie des Sciences de Paris, ser. Mathematiques, 
   submitted in 1995 

\bibitem{d4p} Zapatrin, R.R., 
{\it Algebraic duality in the theory of partially ordered sets},
Order, submitted in 1996 

\end{thebibliography}
\end{document}